\documentclass[letterpaper, 10pt, journal, twocolumn, final]{IEEEtran}

\usepackage{amsthm}
\usepackage{soul}
\usepackage[normalem]{ulem}
\usepackage[hidelinks]{hyperref}
\usepackage{bookmark}
\usepackage[latin1]{inputenc}
\usepackage{microtype}
\usepackage[noadjust]{cite}

\usepackage[dvipsnames]{xcolor}
\usepackage{amsfonts,amssymb,amsmath,color,amsthm}
\usepackage{placeins}
\usepackage{graphicx}
\usepackage{algorithm}
\usepackage[]{algpseudocode}

\usepackage{tikz}
\usetikzlibrary{shapes,arrows}

\usepackage{array}
\usepackage{booktabs}

\usepackage{dsfont} %

\newcommand{\map}[3]{#1: #2 \rightarrow #3}

\newcommand{\EE}{\mathcal{E}} 
\newcommand{\GG}{\mathcal{G}}

\newcommand{\YY}{\mathcal{Y}}

\newcommand{\lexmin}{\mathop{\rm lex\text{-}min}}

\newcommand{\subj}{\textnormal{subj. to}}

\newcommand{\nbrs}{\mathcal{N}}

\newcommand{\real}{{\mathbb{R}}}

\newcommand{\integer}{{\mathbb{Z}}}
\newcommand{\conv}[1]{\textnormal{conv}(#1)}

\newcommand{\until}[1]{\{1,\ldots,#1\}}

\newcommand\oprocendsymbol{\hbox{$\square$}}
\newcommand\oprocend{\relax\ifmmode\else\unskip\hfill\fi\oprocendsymbol}
\def\eqoprocend{\tag*{$\square$}}

\newtheorem{theorem}{Theorem}[section]
\newtheorem{proposition}[theorem]{Proposition}
\newtheorem{corollary}[theorem]{Corollary}
 \newtheorem{lemma}[theorem]{Lemma}
\newtheorem{remark}[theorem]{Remark}

\newtheorem{assumption}[theorem]{Assumption}

\newcommand{\0}{0}
\newcommand{\1}{\mathds{1}}

\newcommand{\bx}{x}
\newcommand{\bz}{z}

\newcommand{\by}{y}

\newcommand{\bmu}{\mu}

\newcommand{\blambda}{\lambda}

\newcommand{\smallsum}{\textstyle\sum\limits}

\newcommand{\tx}{\bx^{\textsc{L}}}

\newcommand{\algx}{\bx^\infty}
\newcommand{\algrho}{\rho^\infty}
\newcommand{\hz}{\widehat{\bz}}
\renewcommand{\restriction}{\sigma}
\newcommand{\asyrestriction}{\restriction^\infty}
\newcommand{\ftrestriction}{\restriction^\textsc{ft}}
\newcommand{\infeas}{\textsc{infeas}}
\newcommand{\feas}{\textsc{feas}}
\renewcommand{\Xi}{X_i^{\textsc{milp}}}

\newcommand{\dz}{p}
\newcommand{\dr}{q}

\newcommand{\INT}{I_{\mathbb{Z}}}
\newcommand{\NINT}{I_{\mathbb{R}}}

\newcommand{\LP}{\textsc{lp}}
\newcommand{\slater}{\textsc{sl}}
\newcommand{\MILP}{\textsc{milp}}

\newcommand{\bL}{L}

\makeatletter
\newcommand{\StatexIndent}[1][3]{%
  \setlength\@tempdima{\algorithmicindent}%
  \Statex\hskip\dimexpr#1\@tempdima\relax}
\makeatother

\renewcommand{\lim}{\operatornamewithlimits{lim\vphantom{p}}}

\graphicspath{{figs/},{figs/simulations/}}

\def \algname/{Distributed Primal Decomposition for MILPs}

\begin{document}

\title{\huge Distributed Primal Decomposition for Large-Scale MILPs}

\author{Andrea Camisa, \IEEEmembership{IEEE Student Member}, 
  Ivano Notarnicola, \IEEEmembership{IEEE Member},
  Giuseppe Notarstefano, \IEEEmembership{IEEE Member} \vspace*{-1cm}
  \thanks{
  A preliminary short version of this paper has appeared as~\cite{camisa2018primal}.
  The current article provides a more detailed discussion, %
  an improved algorithm with tighter restriction, all the theoretical proofs
  and an extensive numerical study.
  }
  \thanks{A. Camisa, I. Notarnicola and G. Notarstefano are with the Department of Electrical, 
  Electronic and Information Engineering, University of Bologna, Bologna, Italy. 
  \texttt{\{a.camisa, ivano.notarnicola, giuseppe.notarstefano\}@unibo.it}.
  This result is part of a project that has received funding from the European 
  Research Council (ERC) under the European Union's Horizon 2020 research 
  and innovation programme (grant agreement No 638992 - OPT4SMART).
  }
}

\maketitle

\begin{abstract}
  This paper deals with a distributed Mixed-Integer Linear Programming (MILP)
  set-up arising in several control applications.
  Agents of a network aim to minimize the sum of local linear cost functions subject
  to both individual constraints and a linear coupling constraint involving all the decision
  variables. A key, challenging feature of the considered set-up is that some components of
  the decision variables must assume integer values.
  The addressed MILPs are NP-hard, nonconvex and large-scale. Moreover, several
  additional challenges arise in a distributed framework due to the coupling constraint,
  so that feasible solutions with guaranteed suboptimality bounds are of interest.
  We propose a fully distributed algorithm based on a primal decomposition approach
  and an appropriate tightening of the coupling constraint.
  The algorithm is guaranteed to provide feasible solutions in finite time.
  Moreover, asymptotic and finite-time suboptimality bounds are established for
  the computed solution.
  Montecarlo simulations highlight the extremely low suboptimality bounds
  achieved by the algorithm.
\end{abstract}


\section{Introduction}
\label{sec:intro}
In this paper, we investigate large-scale Mixed-Integer Linear Programs (MILPs)
that are to be solved by a network of agents without a central coordinator.
The goal is to minimize the sum of objective functions while satisfying individual
constraints and a common, non-sparse coupling constraint.
We term these MILPs \emph{constraint coupled}.
Due to the mixed-integer decision variable, the large-scale size
and the coupling constraint, these problems turn out to be extremely
challenging in a distributed context.
This typically arise in several relevant control applications, such as
microgrid control, economic dispatch in power systems, task assignment
in cooperative robotics.
An interesting scenario arises in distributed Model Predictive Control (MPC),
where a large set of nonlinear control systems must cooperatively solve a
common control task and their states, outputs and/or inputs are coupled through
a coupling constraint. Here, the constraint-coupled structure
results directly from the problem formulation, and the integrality constraints
stem from the MILP approximation of the original optimal control problem
\cite{bemporad1999control}.
In cooperative MPC schemes, such a complex optimization problem should
be ideally solved at each control iteration (see e.g.
\cite{kuwata2010cooperative,
muller2011general,
dinh2013dual,
giselsson2013accelerated,
wang2018accelerated}).
Being these problems NP-hard, it is not computationally affordable to achieve
exact optimality, however feasible (suboptimal) solutions are often sufficient
to guarantee stability.
It is thus of great interest to compute ``good-quality'' feasible solutions
of large-scale MILPs.

Since our paper deals with constraint-coupled optimization,
we organize the relevant literature in two parts.
First, we review existing methods for convex constraint-coupled problems.
In the tutorial paper~\cite{%
necoara2011parallel},
parallel decomposition techniques are reviewed.
A distributed gradient descent method is proposed
in~\cite{lakshmanan2008decentralized} to solve smooth resource allocation
problems.
In~\cite{simonetto2012regularized} a regularized saddle-point algorithm for
convex optimization problems over networks is analyzed.
In \cite{cherukuri2015distributed,cherukuri2016initialization} distributed algorithms based on
Laplacian-gradient dynamics are used to solve economic dispatch over
digraphs.
In~\cite{simonetto2016primal,falsone2017dual} distributed dual decomposition-based algorithms
for constraint-coupled problems are analyzed, while \cite{notarnicola2017constraint}
proposes a distributed algorithm based on successive duality steps.
Approaches for constraint coupled problems based on augmented
Lagrangian methods with consensus schemes are investigated in
\cite{zhang2018consensus,falsone2019tracking}.
Finally, a discussion on approaches for convex constraint-coupled problems can
be found in the recent survey paper~\cite{notarstefano2019distributed}.

Second, we review parallel and distributed algorithms for
MILPs. In~\cite{kim2013scalable} a Lagrange relaxation approach is applied to
demand response control in smart grids.
In~\cite{takapoui2017simple} a heuristic for embedded mixed-integer programming
is studied to obtain approximate solutions.
A first attempt to obtain a distributed approximate solution for MILPs
is~\cite{kuwata2011cooperative}.
Recently, a distributed algorithm %
has been investigated in~\cite{testa2019distributed} to solve a different class of
MILPs with shared decision variable.
A pioneering work on fast, master-client parallel algorithms to
find approximate solutions of our problem set-up
is~\cite{vujanic2016decomposition}.
In~\cite{falsone2019decentralized}, an enhanced
version has been proposed to improve
the quality of the solution. A distributed implementation
of~\cite{falsone2019decentralized} is proposed in~\cite{falsone2018distributed}.

The contributions of this paper are as follows. We first focus on
constraint-coupled convex problems and provide a distributed
algorithm based on a combined relaxation and primal decomposition approach.
Thanks to this first analysis, as second and main contribution of
the paper, we then propose a distributed optimization algorithm for the fast
computation of feasible solutions to
large-scale, constraint-coupled MILPs.
Our algorithm builds on the distributed primal decomposition applied to
a convex approximation of the target MILP with restricted coupling constraint.
For the mixed-integer solution estimate computed by the proposed distributed method,
we are able to: (i) establish both
asymptotic and finite-time feasibility,
and (ii) provide both asymptotic and finite-time suboptimality bounds.
Thanks to the primal decomposition reformulation, the proposed restriction of
the coupling constraint turns out to be tighter than the state of the art.
Through an extensive numerical study on randomly
generated MILPs, we show that our approach is able to achieve interestingly low
suboptimality gaps.

The paper unfolds as follows.
In Section~\ref{sec:set-up_and_preliminaries}, we introduce the MILP
set-up together with useful preliminaries.
In Section~\ref{sec:algorithm_milp}, we propose our distributed algorithm
which is analyzed in Section~\ref{sec:analysis}.
In Section~\ref{sec:simulations}, a numerical study is presented.
All the proofs of the theoretical results are deferred to the appendix.

\section{Optimization Set-up and Preliminaries}
\label{sec:set-up_and_preliminaries}

In this section, we introduce the problem set-up together with some preliminaries that 
act as building blocks for the development of our methodology.

\subsection{Constraint-Coupled Mixed-Integer Linear Program}
\label{sec:set-up}

Let us consider a network of $N$ agents aiming to solve the
optimization problem
\begin{align}
\label{eq:MILP}
\begin{split}
  \min_{\bx_1,\ldots,\bx_N} \: & \: \smallsum_{i =1}^N c_i^\top \bx_i
  \\
  \subj \: 
  & \: \smallsum_{i=1}^N A_i \bx_i \leq b
  \\
  & \: \bx_i \in \Xi, \hspace{1.5cm} i \in \until{N},
\end{split}
\end{align}
where, for all $i \in \until{N}$, the decision variable $\bx_i$ has $\dz_i + \dr_i$
components and the mixed-integer constraint set %
is of the form
$\Xi = P_i \cap (\integer^{\dz_i} \times \real^{\dr_i})$,
for some nonempty compact polyhedron $P_i \subset \real^{\dz_i + \dr_i}$.
The decision variables are intertwined by $S$ linear \emph{coupling} constraints,
described by the matrices $A_i \in \real^{S \times (\dz_i+\dr_i)}$ %
and the vector $b \in \real^S$.
We assume that problem~\eqref{eq:MILP} is feasible
and denote by $(\bx_1^\star, \ldots, \bx_N^\star)$ an optimal solution
with cost $J^\MILP = \sum_{i=1}^N c_i^\top \bx_i^\star$.
In many control applications, %
the number of decision variables is typically much larger than the number of coupling
constraints. Therefore, in this paper we let $N \gg S$, leading to
\emph{large-scale} instances of problem~\eqref{eq:MILP}.

We assume each agent $i$ has a \emph{partial knowledge} of problem~\eqref{eq:MILP},
i.e., it knows only its local data $\Xi$, $c_i$, $A_i$ and $b$.
The goal for each agent is to compute its portion $\bx_i^\star$ of an optimal
solution of~\eqref{eq:MILP}, by means of neighboring communication
and local computation.
Agents communicate according to a connected and undirected graph
$\GG = (\until{N}, \EE)$, where $\EE\subseteq \until{N}^2$ is the
set of edges. 
The set of neighbors of $i$ in $\GG$ is 
$\nbrs_i = \left\{j \in \until{N} \mid (i,j) \in \EE \right\}$.

To solve MILP~\eqref{eq:MILP}, one can employ enumeration
schemes, such as branch-and-bound or cutting-plane techniques.
However, this would not exploit its separable structure and
would result into a computationally intensive algorithm.
Therefore, in the next subsection
we introduce an approximate version of the problem that preserves its structure
while allowing for the application of efficient decomposition techniques.

\subsection{Linear Programming Approximation of the Target MILP}%
\label{sec:lp_approximation}
Following~\cite{vujanic2016decomposition,falsone2018distributed,falsone2019decentralized,bertsekas1982constrained},
let us consider a modified version of problem~\eqref{eq:MILP}
where the right-hand side of the coupling constraint is
restricted by a vector $\restriction \ge \0$ and each mixed-integer
set $\Xi$ is replaced by its convex hull, denoted by $\conv{\Xi}$.
The following convex problem is obtained
\begin{align}
\label{eq:LP_restricted}
\begin{split}
  \min_{\bz_1,\ldots,\bz_N} \: & \: \smallsum_{i =1}^N c_i^\top \bz_i
  \\
  \subj \: 
  & \: \smallsum_{i=1}^N A_i \bz_i \leq b - \restriction
  \\
  & \: \bz_i \in \conv{\Xi}, \hspace{0.7cm} i \in \until{N},
\end{split}
\end{align}
where $\bz_i \in \real^{\dz_i + \dr_i}$ for all $i \in \until{N}$. We introduced
$\bz_i$ to clearly distinguish continuous variables from their mixed-integer
counterparts in problem~\eqref{eq:MILP}.
When $\restriction = \0$, problem~\eqref{eq:LP_restricted} is
a relaxation of problem~\eqref{eq:MILP} and preserves its feasibility.
When $\restriction > 0$, as done in related approaches
\cite{vujanic2016decomposition,falsone2018distributed,falsone2019decentralized},
feasibility of problem~\eqref{eq:LP_restricted} must be also assumed (see
Assumption~\ref{ass:uniqueness}).

The main point in
solving the (convex) problem~\eqref{eq:LP_restricted}
in place of the (nonconvex) original MILP~\eqref{eq:MILP} is to reconstruct
a feasible solution of~\eqref{eq:MILP} starting from a solution of~\eqref{eq:LP_restricted}.
The restriction $\sigma$ is designed to guarantee that the
solution is feasible for the coupling constraint
and depends on the specific mixed-integer reconstruction procedure.
Due to the feasibility assumption, the larger is $\restriction$,
the narrower is the class of problems to which the approach can be applied.
Notably, our method is no worse than \cite{vujanic2016decomposition} in terms of
needed $\restriction$. In fact, numerical experiments highlight an extremely lower
restriction magnitude,
which, as a byproduct, results also in much less suboptimality
of the computed solution.

Next, we introduce a key result based on the Shapley-Folkman lemma.
\begin{proposition}
\label{prop:LP_integer_components}
	Let problem~\eqref{eq:LP_restricted} be feasible and let
	$(\bar{\bz}_1, \ldots, \bar{\bz}_N)$	be any vertex of its feasible set.
	Then, there exists an index set $\INT \subseteq \until{N}$,
	with cardinality $|\INT| \ge N-S$,
	such that $\bar{\bz}_i \in \Xi$ for all $i \in \INT$.
	\oprocend
\end{proposition}
A proof is given in the early reference~\cite{bertsekas1982constrained}
(see also~\cite{vujanic2016decomposition} for a more recent proof).
Since Proposition~\ref{prop:LP_integer_components} gives a bound on the
number of portions that are not mixed integer,
the leading idea is to adapt only these portions while keeping the others intact.
The approach will heavily rely
on a primal decomposition framework whose foundations are reviewed in the next
subsection.

\subsection{Primal Decomposition Review}
\label{sec:primal_decomposition}

Primal decomposition is a powerful tool to recast constraint-coupled convex
programs such as~\eqref{eq:LP_restricted} into a master-subproblem architecture \cite{silverman1972primal}.
The right-hand side vector $b - \restriction$ of the coupling constraint
is interpreted as a given (limited) resource to be shared among agents.
Thus, local \emph{allocation vectors} $\by_i \in \real^S$ for all $i$
are introduced such that $\sum_{i=1}^N \by_i = b - \restriction$. %
To determine the allocations, a \emph{master problem} is introduced
\begin{align}
\begin{split}
  \min_{\by_1,\ldots,\by_N} \: & \: \smallsum_{i =1}^N p_i (\by_i) 
  \\
  \subj \: & \: \smallsum_{i=1}^N \by_i = b - \restriction
  \\
  & \: \by_i \in Y_i, \hspace{1.2cm} i \in\until{N},
\end{split}
\label{eq:primal_decomp_master}
\end{align}
where, for each $i\in\until{N}$, $\map{p_i}{\real^S}{\real}$ is
defined as the optimal cost of the $i$-th (linear programming) \emph{subproblem}
\begin{align}
\begin{split}
  p_i(\by_i) = \: \min_{\bz_i} \: & \: c_i^\top \bz_i
  \\
  \subj \: 
  & \: A_i \bz_i \leq \by_i
  \\
  & \: \bz_i \in \conv{\Xi}.
\end{split}
\label{eq:primal_decomp_subproblem}
\end{align}
In problem~\eqref{eq:primal_decomp_master}, the constraint
$Y_i \subseteq\real^S$ is
the set of $\by_i$ for which
problem~\eqref{eq:primal_decomp_subproblem} is feasible, i.e.,
such that there exists $\bz_i \in \conv{\Xi}$ satisfying the local
\emph{allocation constraint} $A_i \bz_i \le \by_i$.
The following lemma (\cite[Lemma 1]{silverman1972primal}) establishes the equivalence between
LP~\eqref{eq:LP_restricted} and
problems~\eqref{eq:primal_decomp_master}--\eqref{eq:primal_decomp_subproblem}.
\begin{lemma}
\label{lemma:primal_decomp}
  Let~\eqref{eq:LP_restricted} be feasible.
  Then,
    (i) the optimal costs of problems~\eqref{eq:LP_restricted} and~\eqref{eq:primal_decomp_master} are equal;
    (ii) if $(\by_1^\star, \ldots, \by_N^\star)$ is optimal for~\eqref{eq:primal_decomp_master}
    and, for all $i$,
		  $\bz_i^\star$ is optimal for~\eqref{eq:primal_decomp_subproblem}
		  (with $\by_i = \by_i^\star$), then $(\bz_1^\star, \ldots, \bz_N^\star)$
		  is an optimal solution of~\eqref{eq:LP_restricted}.
    \oprocend
\end{lemma}

Note that, given an optimal allocation $(\by_1^\star, \ldots, \by_N^\star)$,
each node can retrieve its portion of an optimal solution of
problem~\eqref{eq:LP_restricted} by relying only on its local allocation
$\by_i^\star$.

\section{\algname/}
\label{sec:algorithm_milp}

In this section we propose a novel distributed algorithm to compute a feasible
solution of MILP~\eqref{eq:MILP}. A cornerstone of the proposed strategy
is the distributed primal decomposition method to solve the convex
problem~\eqref{eq:LP_restricted}. We first present this scheme.
Then, we formally describe our algorithm for MILPs together with its underlying rationale.

\vspace{-0.2cm}
\subsection{Distributed Primal Decomposition for Convex Problems}
\label{sec:algorithm_convex-distributed}
As already discussed, agents can compute
the solution to~\eqref{eq:LP_restricted} by independently solving
problem~\eqref{eq:primal_decomp_subproblem}, provided that
an optimal allocation $(\by_1^\star, \ldots, \by_N^\star)$
is available.
To compute such optimal allocation, one can think of applying a
projected subgradient method to
problem~\eqref{eq:primal_decomp_master}.
By denoting $p(\by) = \sum_{i=1}^N p_i(\by_i)$ the cost function of
problem~\eqref{eq:primal_decomp_master} and by letting $t \ge 0$
be an iteration index, the update reads
$\by^{t+1} = [\by^t - \alpha^t \nabla p(\by^t)]_\YY$, where $\alpha^t$ is the
step-size and $[\cdot]_\YY$ denotes the Euclidean projection onto 
the feasible set of problem~\eqref{eq:primal_decomp_master}.
The constraints $y_i \in Y_i$ in $\YY$ can be handled by resorting 
to a relaxation approach similar to the one used in~\cite{notarnicola2017constraint}. 
As for the constraint $\sum_{i=1}^N y_i = b - \restriction$, the projection
admits a closed-form solution that reads
\begin{align}
  \by_i^{t+1}
  = \by_i^t
  - \frac{\alpha^t}{N} \smallsum_{j=1}^N \big( \nabla p_i(\by_i^t) - \nabla p_j(\by_j^t) \big),
  \hspace{0.8cm} \forall \: i
\label{eq:centr_subgr_method}
\end{align}
where $\alpha^t$ is the step-size.
However, the update~\eqref{eq:centr_subgr_method}
requires knowledge of the entire vector $(\nabla p_1(\by_1^t), \ldots, \nabla p_N(\by_N^t))$
and therefore cannot be performed in a distributed way. We next present an
effective approach that bridges this gap.

Each agent $i$ maintains a local allocation estimate $\by_i^t$.
At each iteration $t \ge 0$, agents compute $\bmu_i^t$ as a Lagrange
multiplier of
\begin{align} 
\label{eq:dpd_z_LP}
\begin{split}
	\min_{\bz_i, v_i} \hspace{1.2cm} &\: c_i^\top \bz_{i} + M v_i
	\\
	\subj \hspace{0.3cm} 
	\: \bmu_i : \: & \: A_i \bz_i \leq \by_i^t + v_i\1
	\\
	& \: \bz_i \in \conv{\Xi}, \: \: v_i \ge 0,
\end{split}
\end{align}
where $M > 0$ and $\1$ is the vector of ones. Then, each agent $i$ receives $\bmu_j^t$ from its
neighbors $j \in \nbrs_i$ and updates $\by_i^t$ with
\begin{align} 
\label{eq:dpd_y_update}
  \by_i^{t+1} = \by_i^t + \alpha^t \smallsum_{j \in \nbrs_i} \big( \bmu_i^t - \bmu_j^t \big),
\end{align}
where $\alpha^t$ is the step-size.
Note that problem~\eqref{eq:dpd_z_LP} is a modified version of~\eqref{eq:primal_decomp_subproblem}
that is feasible for all $y_i$ and that $\nabla p_i(\by_i^t) = \mu_i^t$ \cite[Section 5.4.4]{bertsekas1999nonlinear}.
It is introduced to take into account the 
constraints $y_i \in Y_i$ in~\eqref{eq:primal_decomp_master}. %
Intuitively, update~\eqref{eq:dpd_y_update} is obtained by making the
centralized update~\eqref{eq:centr_subgr_method} match the graph sparsity.
As for the step-size in~\eqref{eq:dpd_y_update}, the following standard assumption is made.
\begin{assumption}
  \label{ass:step-size}
  The step-size sequence $\{ \alpha^t \}_{t\ge0}$, with each $\alpha^t \ge 0$,
  satisfies $\sum_{t=0}^{\infty} \alpha^t = \infty$,
  $\sum_{t=0}^{\infty} \big( \alpha^t \big)^2 < \infty$.
  \oprocend
\end{assumption}

Next, we formalize the convergence result properties
of the distributed primal decomposition algorithm for convex problems
described by~\eqref{eq:dpd_z_LP}--\eqref{eq:dpd_y_update}.
\begin{proposition}
\label{prop:DPD_convergence}
	Let Assumption~\ref{ass:step-size} hold. Moreover, let
	problem~\eqref{eq:LP_restricted} be feasible and let the local allocation
	vectors $\by_i^0$ be initialized such that
	$\sum_{i=1}^N \by_i^0 = b - \restriction$.
	Then, for a sufficiently large $M > 0$, the distributed
	algorithm~\eqref{eq:dpd_z_LP}--\eqref{eq:dpd_y_update} generates an
	allocation vector sequence
	$\{ \by_1^t,\ldots,\by_N^t \}_{t\ge 0}$ and a primal sequence
	$\{ \bz_1^t,\ldots,\bz_N^t \}_{t\ge 0}$ (solutions of problem~\eqref{eq:dpd_z_LP} for all $t \ge0$) 
	such that
  \begin{enumerate}
	  \item $\sum_{i=1}^N \by_i^t  = b - \restriction$, for all $t \ge 0$;
	   
	  \item $\lim_{t \to \infty} \| \by_i^t - \by_i^\star \| = 0$ for all
	    $i \in \until{N}$, where $(\by_1^\star, \ldots, \by_N^\star)$ is an optimal
	    solution of~\eqref{eq:primal_decomp_master};
	  
	  \item every limit point of $\{ \bz_1^t,\ldots,\bz_N^t \}_{t\ge 0}$
	    is an optimal (feasible) solution of problem~\eqref{eq:LP_restricted}.
	    \oprocend
  \end{enumerate}
\end{proposition}
The proof relies on a primal decomposition reinterpretation of the
algorithm in~\cite{notarnicola2017constraint} and is omitted for space constraints.
Operatively, the parameter $M > 0$ can be chosen by using an iterative update scheme as
proposed in~\cite{bertsekas1982constrained}.

\vspace{-0.3cm}
\subsection{Distributed Algorithm Description}
\label{sec:algorithm_milp_description}

We are now ready to formally introduce our \algname/.
Each agent $i$ maintains a mixed-integer solution estimate $\bx_i^t \in \Xi$
and a local allocation estimate $\by_i^t \in \real^S$.
At each iteration $t \ge 0$, the local allocation estimate $\by_i^t$
is updated according to~\eqref{eq:dpd_z_LP}--\eqref{eq:dpd_y_update}.
After $T_f > 0$ iterations, agents compute
a tentative mixed-integer solution $\bx_i^{T_f}$, based
on the last computed allocation $\by_i^{T_f}$ (cf.~\eqref{eq:alg_lexmin_MILP}).
Here the notation $\lexmin$ denotes that $\rho_i, \xi_i$ and $\bx_i$ are
minimized in a lexicographic order \cite{notarstefano2019distributed}. In Section~\ref{sec:mixed_integer_sol_computation},
we discuss in more detail the meaning of problem~\eqref{eq:alg_lexmin_MILP}
and an operative way to solve it.
Algorithm~\ref{alg:algorithm} summarizes the steps from the
perspective of agent $i$.
\begin{algorithm}[H]
\floatname{algorithm}{Algorithm}

  \begin{algorithmic}[0]

    \Statex \textbf{Initialization}:
      Set $T_f > 0$,
      $\by_{i}^0$ such that $\sum_{i=1}^N \by_i^0 = b - \restriction$
    \medskip

    \Statex \textbf{Evolution}: For $t = 0, 1, \ldots, T_f$ \textbf{do} %
    \smallskip
    
      \StatexIndent[0.75]
      \textbf{Compute} $\bmu_i^t$ as a Lagrange multiplier of~\eqref{eq:dpd_z_LP}
      \smallskip

      \StatexIndent[0.75]
      \textbf{Receive} $\bmu_{j}^t$ from $j\in\nbrs_i$ and update $\by_{i}^{t+1}$ with~\eqref{eq:dpd_y_update}
      \medskip
      
      \Statex
      \textbf{Return} $\bx_i^{T_f}$ as optimal solution of
      \begin{align}
      \label{eq:alg_lexmin_MILP}
      \begin{split}      
        \lexmin_{\rho_i, \xi_i, \bx_i} \: & \: \rho_i
        \\
        \subj \: 
        & \: c_i^\top \bx_i \leq \xi_i
        \\
        & \: A_i \bx_i \leq  \by_i^{T_f} + \rho_i \1
        \\
        & \: \bx_i \in \Xi, \: \: \rho_i \ge 0.
      \end{split}
      \end{align}
      \vspace{-0.3cm}

  \end{algorithmic}
  \caption{\algname/}
  \label{alg:algorithm}
\end{algorithm}

A sensible choice for $\by_i^0$ is $\by_i^0 = (b-\restriction)/N$.
By Proposition~\ref{prop:DPD_convergence}, the local allocation
vectors $\{\by_1^t, \ldots, \by_N^t\}_{t \ge 0}$ converge asymptotically to an
optimal solution $(\by_1^\star, \ldots, \by_N^\star)$ of
problem~\eqref{eq:primal_decomp_master}.
Moreover, owing to Proposition~\ref{prop:LP_integer_components}, the asymptotic
solution $\bz_i^\star$ of problem~\eqref{eq:dpd_z_LP} is already
mixed integer for at least $N-S$ agents.
As for the remaining (at most) $S$ agents, a recovery procedure
is needed to guarantee that they also have have a mixed-integer
solution. This is done via step~\eqref{eq:alg_lexmin_MILP}.
In order to allow for a premature halt of the algorithm, we let
\emph{all the agents} perform~\eqref{eq:alg_lexmin_MILP}.
In Section~\ref{sec:analysis}, an asymptotic and finite-time analysis
of Algorithm~\ref{alg:algorithm} is provided.
From an implementation point of view, an explicit description of
$\conv{\Xi}$ in terms of inequalities might not be available.
Nevertheless, agents can still obtain an estimate of $\bmu_i^t$ by
\emph{locally} running a dual subgradient method
on~\eqref{eq:dpd_z_LP}, which involves the solution of
small (local) MILPs without needing $\conv{\Xi}$.
The main limitation of the algorithm is only due
to the local computation capability of each node.

\vspace{-0.2cm}
\subsection{Discussion on Mixed-Integer Solution Recovery}
\label{sec:mixed_integer_sol_computation}

In this subsection, we describe in more detail the %
approach
behind problem~\eqref{eq:alg_lexmin_MILP} and how it allows agents to
recover a ``good'' solution of MILP~\eqref{eq:MILP}.
We first describe the procedure at steady state and then show
how we cope with the finite number of iterations $T_f$.

\begin{figure}[htbp]\centering\vspace{-0.25cm}
  \includegraphics{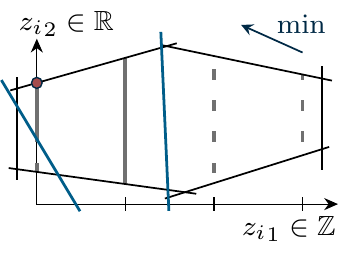}
  \hfill
  \includegraphics{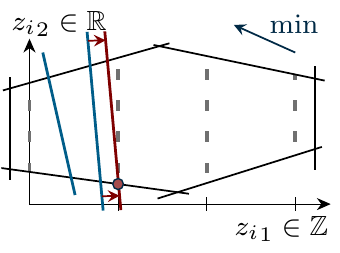}
  
  \vspace{-0.45cm}
  \caption{
    Black lines represent $P_i$ while the blue lines are the two components of 
    $A_i \bz_i \le \by_i^\star$. The resulting feasible set is denoted with solid vertical lines.
    Left: the problem is feasible and admits an %
    optimal solution (red dot), which is mixed integer.
    Right: the problem is infeasible and
    the constraint $A_i \bz_i \le \by_i^\star$ is enlarged just enough to include the closest
    mixed-integer vector (red dot).
  }
\label{fig:solution_computation}
\end{figure}

Let $(\bz_1^\star, \ldots, \bz_N^\star)$ be an optimal solution of the approximate
problem~\eqref{eq:LP_restricted} and let $(\by_1^\star, \ldots, \by_N^\star)$ be
a corresponding allocation of the master
problem~\eqref{eq:primal_decomp_master}, computed asymptotically by
Algorithm~\ref{alg:algorithm}.
A straight approach to recover a mixed-integer solution would be to
solve for all $i$ the optimization problem
\begin{align}
\begin{split}
  \min_{\bx_i} \: & \: c_i^\top \bx_i
  \\
  \subj \: 
  & \: A_i \bx_i \leq \by_i^\star
  \\
  & \: \bx_i \in \Xi.
\end{split}
\label{eq:local_MILP_description}
\end{align}
Problem~\eqref{eq:local_MILP_description} is a (small) local MILP that
is reminiscent of the subproblem~\eqref{eq:primal_decomp_subproblem}.
Depending on the allocation constraint $A_i \bx_i \le \by_i^\star$,
problem~\eqref{eq:local_MILP_description} may be feasible or not.
Figure~\ref{fig:solution_computation} (left) shows an example
with $z_i \in \integer \times \real$, in which~\eqref{eq:local_MILP_description} is feasible.

\noindent%
In view of Proposition~\ref{prop:LP_integer_components},
at least $N-S$ portions $z_i^\star$ of optimal solution
are already mixed-integer and thus optimal for the corresponding local
MILPs~\eqref{eq:local_MILP_description}.
Recall that $N \gg S$, thus \emph{the majority} of
subproblems~\eqref{eq:local_MILP_description} are feasible,
while the total number of
(possibly) infeasible instances is \emph{at most} $S$.
If~\eqref{eq:local_MILP_description} is infeasible for some agent $i$, we let
the agent find a mixed-integer vector with the
\emph{minimal} violation
of the allocation constraint as depicted in
Figure~\ref{fig:solution_computation} (right).

The procedure just outlined for the asymptotic case %
is entirely encoded in problem~\eqref{eq:alg_lexmin_MILP}.
To see this, let us now show how to operatively solve~\eqref{eq:alg_lexmin_MILP}
with $\by_i^\star$ in place of $\by_i^{T_f}$.
First, the needed violation of the allocation constraint is determined by
computing $\algrho_i$ as the optimal cost of
\begin{align}
\label{eq:rho_MILP_description}
\begin{split}
  \min_{\rho_i, \bx_i} \: & \: \rho_i
  \\
  \subj \: 
  & \: A_i \bx_i \leq  \by_i^\star + \rho_i \1
  \\
  & \: \bx_i \in \Xi, \: \: \rho_i \ge 0.
\end{split}
\end{align}
Then, the value of $\rho_i$ is fixed to $\algrho_i$ and problem~\eqref{eq:rho_MILP_description}
is re-optimized with its cost function %
replaced by $c_i^\top \bx_i$ to compute $\algx_i$.
If problem~\eqref{eq:local_MILP_description} is feasible,
then %
$\algrho_i = 0$ and the procedure
is equivalent to solving~\eqref{eq:local_MILP_description}.
Instead, if problem~\eqref{eq:local_MILP_description} is infeasible,
a violation $\algrho_i \1$ of the allocation
constraint is permitted. %

Due to the violations, the aggregate mixed-integer solution
  $(\algx_1, \ldots, \algx_N)$ %
  may exceed the \emph{restricted} total resource $b - \restriction$.
  Indeed, although it may not hold $\sum_{i=1}^N A_i \algx_i \le b - \restriction$,
  in the next subsection we show how to design the
  restriction $\restriction$ to ensure that our original goal $\sum_{i=1}^N A_i \algx_i \le b$
  is achieved.

Now, let us discuss how this approach can be adapted to cope with
the finite number of iterations.
The local allocation $\by_i^t$ can be thought of as
$\by_i^t = \by_i^\star + \Delta_i^t$, where $\{\Delta_i^t\}_{t \ge 0} \to 0$ as
$t$ goes to infinity.  By looking at problem~\eqref{eq:rho_MILP_description}, it
is natural to expect that $\rho_i^t \le \algrho_i + \Delta_i^t$. %
Thus, we let \emph{all the agents} perform
step~\eqref{eq:alg_lexmin_MILP} (implemented as in the asymptotic case).
By employing an additional (small) restriction,
we can guarantee that -- after a sufficiently large time -- the total
violation is embedded into the restriction. A detailed analysis of this
approach is given in Section~\ref{sec:finite_time_analysis}.

\vspace{-0.3cm}
\subsection{Design of the Coupling Constraint Restriction}
\label{sec:restriction_approach}

As already discussed, the purpose of the restriction $\sigma$ is to
compensate for possible violations of problem~\eqref{eq:alg_lexmin_MILP}.
Intuitively, we wish to make $\sigma$ as small as possible for two reasons:
the larger is $\restriction$, then (i) the more likely is~\eqref{eq:LP_restricted}
to be infeasible, (ii) the higher is the cost of
the optimal solution $(\bz_1^\star, \ldots, \bz_N^\star)$ of~\eqref{eq:LP_restricted},
which in turn deteriorates the cost of $(\algx_1, \ldots, \algx_N)$.

We now propose a method to find a small a-priori restriction %
to guarantee feasibility of the computed solution.
As before, we focus on the asymptotic case, while the
extension to finite time is given in the following sections.
Intuitively, the restriction must take into account the worst-case
violation due to the mismatch between $(\algx_1, \ldots, \algx_N)$ and
$(\bz_1^\star, \ldots, \bz_N^\star)$.
Such worst case occurs when all the (at most $S$) agents for which
$\bz_i^\star \notin \Xi$ have infeasible instances of~\eqref{eq:local_MILP_description},
leading to a positive violation $\algrho_i > 0$.
Thus, we define the a-priori restriction $\asyrestriction \in \real^S$ as
\begin{align}
  \label{eq:restriction_definition}
  \asyrestriction = S \cdot \max_{i\in\until{N}} \sigma_i^\textsc{loc},
\end{align}
where $\sigma_i^\textsc{loc} \in \real^S$ is the worst-case violation of
agent $i$ and $\max$ is intended component wise (we stick to this convention
from now on).

Let us now quantify $\sigma_i^\textsc{loc}$.
Since $\conv{\Xi}$ is bounded, it is possible to
find a lower-bound vector, which we denote by $\bL_i \in \real^S$, for any 
\emph{admissible} local allocation $\by_i$,
\begin{align}
  \bL_i
  \triangleq \!
  \min_{\bx_i \in \conv{\Xi}} A_i \bx_i
  = \!
  \min_{\bx_i \in \Xi} A_i \bx_i,
\label{eq:minimum_resource_vector}
\end{align}
By construction, it holds
$\bL_i \le A_i \bz_i^\star \le \by_i^\star$.
Recall that each agent computes the needed violation through
problem~\eqref{eq:rho_MILP_description}.
Then, the worst-case violation that may occur at steady-state is
$\rho_i^\textsc{max} \1$, where we define
$\rho_i^\textsc{max}$ as
\begin{align}
\label{eq:minimum_resource_MILP}
\begin{split}
  \rho_i^\textsc{max} \triangleq \min_{\bx_i \in \Xi} \: & \: \max_{s \in \until{S}} \: [A_i \bx_i  - \bL_i]_s,
\end{split}
\end{align}
where the notation $[\cdot]_s$ denotes the $s$-th component of a vector.
Note that the optimization in~\eqref{eq:minimum_resource_MILP} allows each
agent $i$ to find the ``first'' feasible vector, i.e., with minimal resource usage.

\begin{figure}[!htpb]\centering
  \includegraphics[scale=1]{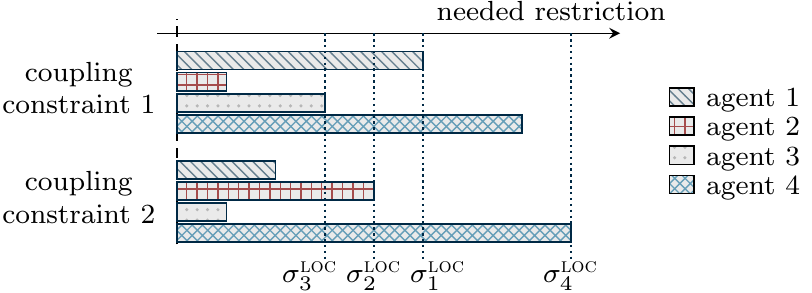}
  \vspace{-0.3cm}
  
  \caption{ Graphical representation of the restriction for $N = 4$ agents and
    $S = 2$ coupling constraints.  Each bar represents the quantity
    $[A_i \bx_i - \bL_i]_s$ at an optimal solution
    of~\eqref{eq:minimum_resource_MILP}, which is the violation associated to
    the first feasible vector with minimal resource usage. The dotted lines
    represent $\rho_i^\textsc{max}$ for all $i$, which here coincide with
    $\restriction_i^\textsc{loc}$.  }
  \label{fig:restriction}
\end{figure}

In order to reduce possible conservativeness of the violation, which can occur
when $\rho_i^\textsc{max} > \max_{\bx_i \in \Xi} [A_i \bx_i - \bL_i]_s$ for some
component $s$ of the coupling constraint, the computation of $\sigma_i^\textsc{loc}$
can be replaced by the saturated version
\begin{align*}
  \sigma_i^\textsc{loc} = \min \Big\{ \rho_i^\textsc{max} \1, \:\: \max_{\bx_i \in \Xi} \: (A_i \bx_i - \bL_i) \Big\}.
\end{align*}
In numerical computations, we have found that usually
$\rho_i^\textsc{max} \ll \max_{\bx_i \in \Xi} \: [A_i \bx_i - \bL_i]_s$,
leading to $\sigma_i^\textsc{loc} = \rho_i^\textsc{max} \1$.
We point out that the computation of $\asyrestriction$ must be performed
in the initialization phase, which can be also carried out in a fully distributed way
by using a max-consensus algorithm.
In Figure~\ref{fig:restriction}, we illustrate an example of the restriction.
\vspace{-0.2cm}

\begin{remark}
\label{rem:comparison_restriction}
  In~\cite{vujanic2016decomposition}, an alternative approach based
  on dual decomposition is explored to compute a feasible solution for
  MILP~\eqref{eq:MILP}.
  The restriction proposed in~\cite{vujanic2016decomposition} is
	\begin{align}
	  \label{eq:vujanic_restriction}
	  \restriction^\textsc{dd} = S \cdot
	  \max_{i \in\until{N}} \: \max_{ \bx_i \in \Xi } (A_i \bx_i - \bL_i),
	\end{align}
	The term $\displaystyle\max_{ \bx_i \in \Xi } (A_i \bx_i - \bL_i)$ may be overly conservative
	and in our approach it
	is replaced with $\sigma_i^\textsc{loc}$,
	which can be thought of as the resource utilization of a feasible vector %
	\emph{closest} to $\bL_i$.
	Independently of the problem at hand, it holds $\asyrestriction \leq \restriction^\textsc{dd}$.
	\oprocend
\end{remark}

\vspace{-0.3cm}
\section{Analysis of Distributed Algorithm}
\label{sec:analysis}
In this section we provide both asymptotic and finite-time analyses of
Algorithm~\ref{alg:algorithm} under the following assumption.
\begin{assumption}%
\label{ass:uniqueness}
  For a given $\restriction \ge 0$, the optimal solution of
  problem~\eqref{eq:LP_restricted} is unique.\oprocend
\end{assumption}%
\noindent
This assumption ensures that the optimal solution of~\eqref{eq:LP_restricted}
is a vertex (hence Proposition~\ref{prop:LP_integer_components} applies).
It can be guaranteed by simply adding a small, random perturbation
to the cost vectors $c_i$.
Notice that Assumption~\ref{ass:uniqueness} is also needed in dual
decomposition approaches such as~\cite{vujanic2016decomposition,falsone2019decentralized,falsone2018distributed}.
Notably, dual decomposition approaches require also uniqueness of the dual
optimal solution of problem~\eqref{eq:LP_restricted}, while our approach is
less restrictive since this is not necessary.

\vspace{-0.35cm}
\subsection{Asymptotic analysis}
\label{sec:asymptotic_analysis}
First, we proceed under the assumption that $\restriction = \asyrestriction$ as
in~\eqref{eq:restriction_definition} and that the algorithm is
executed until convergence to an optimal allocation $(\by_1^\star, \ldots, \by_N^\star)$,
i.e., an optimal solution of problem~\eqref{eq:primal_decomp_master}. Indeed,
steps~\eqref{eq:dpd_z_LP}--\eqref{eq:dpd_y_update} implement
the distributed algorithm in Section~\ref{sec:algorithm_convex-distributed} for
the solution of problem~\eqref{eq:LP_restricted}, so that
Proposition~\ref{prop:DPD_convergence} (ii) applies.

The next theorem shows feasibility of the computed mixed-integer
solution for the target MILP~\eqref{eq:MILP}.
\begin{theorem}[Feasibility]\label{thm:asymptotic_feasibility}
  Let $\restriction = \asyrestriction$ as in~\eqref{eq:restriction_definition},
  and let problem~\eqref{eq:LP_restricted} be feasible and satisfy
  Assumption~\ref{ass:uniqueness}.  Let $(\by_1^\star, \ldots, \by_N^\star)$ be
  an optimal solution of problem~\eqref{eq:primal_decomp_master}.
  Then, the vector $(\algx_1, \ldots, \algx_N)$, with each $\algx_i$
  optimal solution of~\eqref{eq:alg_lexmin_MILP} with
  $\by_i^t = \by_i^\star$, is feasible for MILP~\eqref{eq:MILP}, i.e., $\algx_i \in \Xi$
  for all $i \in \until{N}$ and $\sum_{i=1}^N A_i \algx_i \leq b$.
  \oprocend
\end{theorem}
The proof of Theorem~\ref{thm:asymptotic_feasibility} is given in appendix.

\begin{remark}
  The proof of Theorem~\ref{thm:asymptotic_feasibility} reveals that
	the same result can be obtained by using an allocation
	$(\by_1, \ldots, \by_N)$ associated to \emph{any} vertex of the feasible
	set of problem~\eqref{eq:LP_restricted} (rather than an optimal allocation $(\by_1^\star, \ldots, \by_N^\star)$).
	Proposition~\ref{prop:LP_integer_components}
	can still be applied and the proof remains unchanged.
	\oprocend
\end{remark}

Theorem~\ref{thm:asymptotic_feasibility} guarantees that the computed solution
is feasible for the target MILP~\eqref{eq:MILP}, but, in general,
there is a certain degree of suboptimality.
In the following, we provide suboptimality bounds
under Slater's constraint qualification. %
\begin{assumption}%
\label{ass:slater}
  For a given $\restriction > 0$, there exists a vector $(\hz_1,\ldots, \hz_N)$,
  with $\hz_i \in \conv{\Xi}$ for all $i$, such that
  \begin{align}  
    \zeta & \triangleq \min_{s \in \until{S}} \Big[ b - \restriction -\smallsum_{i=1}^N A_i \hz_i \Big]_s > 0.
    \label{eq:slater_zeta}
  \end{align}
  The cost of $(\hz_1,\ldots, \hz_N)$ is denoted by $J^\slater = \sum_{i=1}^N c_i^\top \hz_i$.
\oprocend
\end{assumption}

\noindent
The following result establishes an a-priori suboptimality bound on the
mixed-integer solution with $\restriction = \asyrestriction$ as
in~\eqref{eq:restriction_definition}.
Due to space constraints, the proofs of
Theorem~\ref{thm:asymptotic_performance_apriori}
and Corollary~\ref{cor:asymptotic_performance} are omitted.
The reader is referred to~\cite{vujanic2016decomposition} for similar results.
\begin{theorem}[A-Priori Suboptimality Bound]
\label{thm:asymptotic_performance_apriori}
  Consider the same assumptions and quantities of
  Theorem~\ref{thm:asymptotic_feasibility}
  and let also Assumption~\ref{ass:slater} hold.
  Then, $(\algx_1, \ldots, \algx_N)$ satisfies the suboptimality bound
  $\sum_{i=1}^N c_i^\top \algx_i - J^\MILP \le B$,
  where $J^\MILP$ is the optimal cost of~\eqref{eq:MILP}
  and $B$ is defined as
  \begin{align*}
    B
    \triangleq
    \bigg( S + \frac{N \| \asyrestriction \|_\infty}{\zeta} \bigg)
      \max_{i \in \until{N}} \: \gamma_i,
  \end{align*}
  with $\zeta$ defined in~\eqref{eq:slater_zeta}, and
  $\displaystyle\gamma_i
    \triangleq \!\!
    \max_{\bx_i \in \Xi} \! c_i^\top \bx_i
    - \!\!\min_{\bx_i \in \Xi} \! c_i^\top \bx_i$.
    \oprocend
\end{theorem}

Note that, although the bound provided by Theorem~\ref{thm:asymptotic_performance_apriori}
is formally analogous to \cite[Theorem 3.3]{vujanic2016decomposition}, there is an implicit
difference due to the restriction values (cf. Remark~\ref{rem:comparison_restriction}).
In particular, our bound is tighter since $\asyrestriction$
is less than or equal to the restriction proposed by~\cite{vujanic2016decomposition}.

A tighter bound can be derived by using the steady-state solution of the
algorithm and computing also the primal solution of~\eqref{eq:dpd_z_LP}
for all $i$. For this reason, we call this bound a posteriori, since it depends on the solution
computed by Algorithm~\ref{alg:algorithm}.
\begin{corollary}[A-Posteriori Suboptimality Bound]
\label{cor:asymptotic_performance}
  Consider the same assumptions and quantities of Theorem~\ref{thm:asymptotic_performance_apriori}.
  Then, %
  $\sum_{i=1}^N c_i^\top \algx_i - J^\MILP \le B^\prime$,
  where $B^\prime$ is defined as
  \begin{align}
    B^\prime
    \triangleq
    \smallsum_{i \in \NINT} ( c_i^\top \algx_i \! - \! c_i^\top \bz_i^\star )
    + \frac{\| \asyrestriction\|_\infty}{\zeta} ( J^\slater \!-\! \smallsum_{i=1}^N \! c_i^\top \bz_i^\star ),
  \label{eq:a_posteriori_bound}
  \end{align}
  with $\zeta$ and $J^\slater$ defined in Assumption~\ref{ass:slater}
  and $\NINT$ containing the indices of agents such that $\bz_i^\star \notin \Xi$
  ($|\NINT| \le S$).
  \oprocend
\end{corollary}

\vspace{-0.3cm}
\subsection{Finite-time Analysis}
\label{sec:finite_time_analysis}

In this section, we provide a finite-time analysis of the distributed algorithm.
All the proofs of this subsection are given in appendix.
To this end, we assume that the restriction is equal to an
enlarged version of the asymptotic restriction in~\eqref{eq:restriction_definition}, i.e.,
\begin{align}
  \label{eq:ftrestriction_definition}
  \ftrestriction = \asyrestriction + \delta\1,
\end{align}
for an arbitrary $\delta > 0$. We assume problem~\eqref{eq:LP_restricted}
is feasible with this new restriction.
We provide two results that extend the results of
Section~\ref{sec:asymptotic_analysis} to a finite-time setting.

At a high level, finite-time feasibility hinges upon the fact
that %
the allocation sequence
$\{\by_1^t, \ldots, \by_N^t\}_{t \ge 0}$ approaches an optimal allocation.
Eventually, the additional restriction $\delta$ can embed the
distance of the current allocation estimate to optimality. The next theorem
formalizes this result.
\begin{theorem}[Finite-time feasibility]
\label{thm:finite_time_feasibility}
  Let $\restriction = \ftrestriction$ as in~\eqref{eq:ftrestriction_definition}, for some
  $\delta > 0$, and let problem~\eqref{eq:LP_restricted} be feasible and satisfy
  Assumption~\ref{ass:uniqueness}.
  Consider
  the mixed-integer sequence $\{\bx_1^t, \ldots, \bx_N^t\}_{t \ge 0}$
  generated by Algorithm~\ref{alg:algorithm} under Assumption~\ref{ass:step-size},
  with $\sum_{i=1}^N \by_i^0 = b - \ftrestriction$.
  There exists a sufficiently large time $T_\delta > 0$ such that the vector
  $(\bx_1^t, \ldots, \bx_N^t)$ is a feasible solution for problem~\eqref{eq:MILP},
  i.e., $\bx_i^t \in \Xi$ for all $i \in \until{N}$ and $\sum_{i=1}^N A_i \bx_i^t \leq b$,
  for all $t \ge T_\delta$.
  \oprocend
\end{theorem}

In principle, the smaller is $\delta$, the longer it takes
for the mixed-integer vector $(\bx_1^t, \ldots, \bx_N^t)$ to satisfy the
coupling constraint. As a function of $\delta$, there is a trade-off between
the number of iterations to guarantee solution feasibility %
and how strict is the assumption that problem~\eqref{eq:LP_restricted} is feasible.

Next, we provide a suboptimality bound that can be evaluated when
the algorithm is halted.
\begin{theorem}[Finite-time suboptimality bound]
\label{thm:finite_time_performance}
  Consider the same assumptions and quantities of Theorem~\ref{thm:finite_time_feasibility}
  and let also Assumption~\ref{ass:slater} hold.
  Moreover, let $\epsilon_i > 0$ for $i \in \until{N}$. %
  Then, there exists a %
  time $T_\epsilon > 0$ such that the vector
  $(\bx_1^t, \ldots, \bx_N^t)$ satisfies the suboptimality bound
  $\sum_{i=1}^N c_i^\top \bx_i^t - J^\MILP \le B^t$ for all $t \ge T_\epsilon$,
  with $B^t$ being
  \begin{align}
    B^t
    \triangleq
    \smallsum_{i=1}^N ( c_i^\top \bx_i^t - J_i^{\LP,t})%
    +
    \smallsum_{i=1}^N \epsilon_i \| \bmu_i^t\|_1 
    +
    \Gamma \|\ftrestriction\|_\infty,
  \label{eq:finite_time_bound}
  \end{align}
  where $J^\MILP$ is the optimal cost of \eqref{eq:MILP},
  $J_i^{\LP,t}$ and $\bmu_i^t$ are the optimal cost and a Lagrange multiplier
  of~\eqref{eq:dpd_z_LP} at time $t$,
  $\Gamma = \frac{N}{\zeta}
    \smallsum_{i=1}^N \Big( \max\limits_{\bx_{i} \in \Xi} c_i^\top \bx_i 
    - \min\limits_{\bx_{i} \in \Xi} c_i^\top  \bx_{i}
    \Big)$, and $\zeta$ associated to any Slater point (cf. Assumption~\ref{ass:slater}).
  \oprocend
\end{theorem}

We point out that, differently from the asymptotic case, the bound~\eqref{eq:finite_time_bound}
does not depend on the solution of~\eqref{eq:dpd_z_LP},
but only on its optimal cost.
Moreover, notice that the bound~\eqref{eq:finite_time_bound} is a posteriori,
while if an a-priori bound with restriction $\ftrestriction$ is desired, it
still has the form of Theorem~\ref{thm:asymptotic_performance_apriori} (since it
does not depend on the algorithmic evolution).

\section{Monte Carlo Numerical Computations}
\label{sec:simulations}

In this section, we provide a computational study on randomly generated MILPs
to compare Algorithm~\ref{alg:algorithm} with~\cite{vujanic2016decomposition}.
The distributed algorithms are emulated using a single machine with the
\textsc{Matlab} software and the local MILPs are solved using the integrated
numerical solver.

We consider large-scale problems with a total of $4500$
optimization variables ($3000$ are integer and
$1500$ are continuous).
There are $N = 300$ agents and $S = 5$
coupling constraints. The local constraints $\Xi$ are subsets of
$\integer^{10} \times \real^5$ satisfying $D_i \bx_i \leq d_i$,
where $D_i$ and $d_i \in \real^{20}$ have random entries in
$[0,1]$ and $[20,40]$ respectively.
Box constraints $-60 \, \1 \le \bx_i \le 60 \, \1$ are added
to ensure compactness.
The cost vector is $c_i = D_i^\top \hat{c}_i$, where
$\hat{c}_i$ has random entries in $[0,5]$.
As for the coupling, the matrices $A_i$ are random with
entries in $[0,1]$ and the resource vector $b \in \real^5$ is random with values in two
different intervals.
Specifically, we first pick values in $[-20N,-15N]$,
which results in a ``loose'' $b$,
then we pick values in $[-180N,-175N]$,
which results in a ``tight'' $b$.

\begin{figure}[htbp]\centering
  \includegraphics[scale=1]{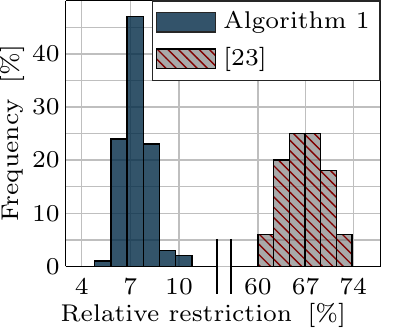}
  \includegraphics[scale=1]{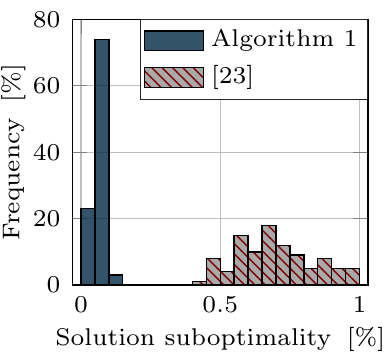}
  
  \caption{
  Montecarlo simulations on random MILPs with loose $b$: comparison histograms
    of Algorithm~\ref{alg:algorithm} and of the method in~\cite{vujanic2016decomposition}.
    See the text for details.}
\label{fig:montecarlo_histogram}
\end{figure}

A total of $100$ MILPs with loose $b$ and $100$ MILPs
with tight $b$ are generated. For each problem,
we check feasibility of problem \eqref{eq:LP_restricted} for both our method
and the method in~\cite{vujanic2016decomposition}.
Then, both algorithms are executed up to asymptotic convergence
to evaluate the mixed-integer solution suboptimality.
The results are summarized in Figures~\ref{fig:montecarlo_table}
and~\ref{fig:montecarlo_histogram}, where the restriction size is
computed as $\|\restriction\| / \|b\|$ and the suboptimality is computed as
$\big(\sum_{i=1}^N c_i^\top \algx_i - q^\star\big)/q^\star$,
with $q^\star$ being the optimal cost of~\eqref{eq:LP_restricted}.
The number of solvable instances is the number of problems for
which problem~\eqref{eq:LP_restricted} is feasible.
For loose $b$, both methods are always applicable.
However, our approach provides better solution
performance than~\cite{vujanic2016decomposition}.
For tight resource vectors, our method is still applicable in the 70\% of
the cases, and provides an average suboptimality of 6.91\%, while the approach
in~\cite{vujanic2016decomposition} cannot be applied due to infeasibility
of the approximate problem~\eqref{eq:LP_restricted} (caused by the
too large needed restrictions).
It is worth noting that the generation interval $[-180N,-175N]$
cannot be further tightened. Indeed, for smaller values of $b$, the target
MILP~\eqref{eq:MILP} becomes infeasible.

\begin{figure}[htbp]\centering
  \includegraphics[scale=1]{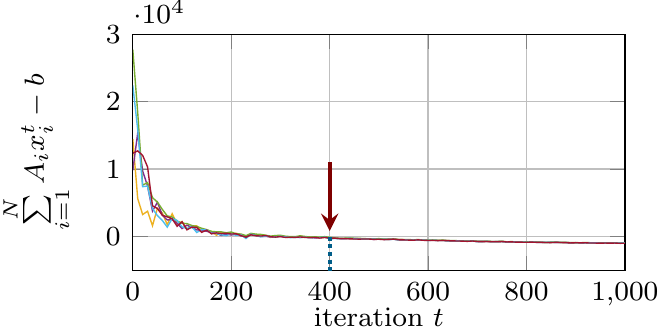}
  \caption{
  Components of the coupling constraint associated to
  $(x_1^t, \ldots, x_N^t)$ generated by Algorithm~\ref{alg:algorithm}
  on a random MILP. The agents are assumed to compute the
  mixed-integer solution according to~\eqref{eq:alg_lexmin_MILP}
  at each iteration.}
\label{fig:alg_evolution}
\end{figure}

Finally, we show the evolution of Algorithm~\ref{alg:algorithm} on a single
instance over an Erd\H{o}s-R\'{e}nyi graph with edge probability $0.2$.
Figure~\ref{fig:alg_evolution} shows the value of the coupling
constraints along the algorithmic evolution, with $\delta = 0.5$
(cf.~\eqref{eq:ftrestriction_definition}). Note that feasibility is achieved
in finite time (within 400 iterations), confirming
Theorem~\ref{thm:finite_time_feasibility}. 
In order to detect feasibility, agents can run
a consensus-based scheme from time to time to check whether the current
solution satisfies the coupling constraints.

\begin{figure}[htbp]\centering\small
  \begin{tabular}{r|cc|cc}
  & \multicolumn{2}{c|}{Algorithm~\ref{alg:algorithm}} & \multicolumn{2}{c}{\cite{vujanic2016decomposition}} \\
  & $b$ loose & $b$ tight & $b$ loose & $b$ tight \\
  \midrule 
  \# solvable problems        & 100\%   & 70\%    & 100\%  & 0\% \\
  size of restriction              & 7.4\%    & 0.72\% & 66.9\% & 6.63\% \\
  suboptimality of solution   & 0.06\%  & 6.91\% & 0.7\%   & N.A. \\ 
  \end{tabular}
  \caption{
    Montecarlo simulations on random MILPs: performance comparison
    of Algorithm~\ref{alg:algorithm} and of the method in~\cite{vujanic2016decomposition}.
    The second and the third row are averaged over the Monte Carlo runs.
    See the text for further details.
  }
\label{fig:montecarlo_table}
\end{figure}

\section{Conclusions}
In this paper we proposed a distributed algorithm to compute a feasible
 solution to large-scale MILPs with high optimality degree. The method is based on
a primal decomposition approach and a suitable restriction of the coupling
constraint and guarantees feasibility of the computed mixed-integer vectors
in finite time.
Asymptotic and finite-time results for feasibility and suboptimality bounds are proved.
Numerical simulations highlight the efficacy of the proposed methodology.

\begin{footnotesize}
\bibliographystyle{IEEEtran}
\bibliography{primal_decomp_milp_biblio}
\end{footnotesize}

\appendix

\section{Proofs}

\subsection{Proof of Theorem~\ref{thm:asymptotic_feasibility}}

  For the sake of analysis, let us denote
  by $(\bz_1^\star, \ldots, \bz_N^\star)$ the optimal solution
  of the restricted LP~\eqref{eq:LP_restricted}.
  By Assumption~\ref{ass:uniqueness}, $(\bz_1^\star, \ldots, \bz_N^\star)$ is a vertex,
  so that by Proposition~\ref{prop:LP_integer_components} %
  there exists $\INT \subseteq\until{N}$, with $|\INT| \ge N-S$, such that
  $\bz_i^\star \in \Xi$ for all $i\in \INT$.
  By Lemma~\ref{lemma:primal_decomp}, $\bz_i^\star$ is an optimal solution of
  problem~\eqref{eq:primal_decomp_subproblem}, with $\by_i = \by_i^\star$,
  for all $i\in\until{N}$.
  Thus, for all $i \in \INT$, $\bz_i^\star \in \Xi$ is the optimal solution
  of~\eqref{eq:alg_lexmin_MILP} with $\by_i^t = \by_i^\star$. Then, it holds
    $A_i \algx_i \le \by_i^\star$ for all $i \in \INT$.

  Let us focus on the set $\NINT = \until{N} \setminus \INT$, which contains indices
  such that $\bz_i^\star \notin \Xi$.
  Let us further partition $\NINT = I_\feas \cup I_\infeas$,
  where the indices collected in $I_\feas$ correspond to feasible
  subproblems~\eqref{eq:primal_decomp_subproblem},
  from which it follows that
$A_i \algx_i \le \by_i^\star$ for all $i \in I_\feas$,
  while the remaining index set $I_\infeas$ corresponds to infeasible
  subproblems~\eqref{eq:primal_decomp_subproblem}.
  We have
  \begin{align*}
    A_i \algx_i 
     \stackrel{(a)}{\le} \by_i^\star + \algrho_i \1
     \stackrel{(b)}{\le} \by_i^\star + \rho_i^\textsc{max} \1,
    \hspace{1cm}
    \forall \: i \in I_\infeas,
  \end{align*}
  where \emph{(a)} follows by construction of $\algx_i$ and \emph{(b)}
  follows since any optimal solution of problem~\eqref{eq:minimum_resource_MILP},
  say $\tx_i$, is feasible for problem~\eqref{eq:rho_MILP_description} (since
  $A_i \tx_i \le \bL_i + \rho_i^\textsc{max} \le \by_i^\star + \rho_i^\textsc{max}$),
  from which it follows that $\algrho_i \le \rho_i^\textsc{max}$ (by optimality).
  Also, notice that, since $\algx_i \in \Xi$, then
  $[A_i \algx_i]_s \le \max_{\bx_i \in \Xi} [A_i \bx_i]_s$ for all $s$ and it holds
  $A_i \algx_i - \by_i^\star
    \le \max_{\bx_i \in \Xi} A_i \bx_i - \by_i^\star
    \le \max_{\bx_i \in \Xi} A_i \bx_i - \bL_i$,
  where $\max$ is intended component wise. It follows that
  $A_i \algx_i  \le \by_i^\star + \sigma_i^\textsc{loc}$ for all $i \in I_\infeas$.
  By summing over $i \in I_\infeas$ the term $\sigma_i^\textsc{loc}$, we obtain
  \begin{align}
    \smallsum_{i \in I_\infeas} \sigma_i^\textsc{loc}
    & 
    \leq |I_\infeas| \: \max_{i \in I_\infeas} \sigma_i^\textsc{loc}
    \nonumber
    \\
    &
    \leq S \max_{i \in \until{N}} \sigma_i^\textsc{loc}
    =
    \asyrestriction,
  \label{eq:i_infeas_sum}
  \end{align}
  where $\max$ is intended component wise.
  Collecting the previous conditions %
  leads to
  \begin{align*}
    \smallsum_{i=1}^N A_i \algx_i
    &=
    \smallsum_{i \in \INT} A_i \algx_i
    + \!\!
    \smallsum_{i \in I_\feas} A_i \algx_i
    + \!\!
    \smallsum_{i \in I_\infeas} A_i \algx_i
    \\
    & 
    \le
    \smallsum_{i=1}^N \by_i^\star
    +
    \smallsum_{i \in I_\infeas} \sigma_i^\textsc{loc}
    \le
    b - \asyrestriction + \asyrestriction = b,
  \end{align*}
  where we used $\sum_{i=1}^N \by_i^\star = b - \asyrestriction$.
  The proof follows.
\oprocend

\subsection{Proof of Theorem~\ref{thm:finite_time_feasibility}}

  Let $\{\by_1^t, \ldots, \by_N^t\}_{t \ge 0}$ denote the allocation vector
  sequence generated by Algorithm~\ref{alg:algorithm}.
  By Proposition~\ref{prop:DPD_convergence}, the sequence
  $\{\by_1^t, \ldots, \by_N^t\}_{t \ge 0}$ converges to an optimal solution
  $(\by_1^\star, \ldots, \by_N^\star)$ of problem~\eqref{eq:primal_decomp_master}
  with $\restriction = \asyrestriction + \delta\1$.
  Thus, for all $i \in \until{N}$ and $\epsilon_i > 0$, there exists
  $T_{\epsilon_i} > 0$ such that $t \ge T_{\epsilon_i} \Rightarrow \|\by_i^t - \by_i^\star\|_\infty \le \epsilon_i$.
  If we let $T = \max_{i \in\until{N}} T_{\epsilon_i}$, then
  $\by_i^t \le \by_i^\star + \epsilon_i\1$ for all $t \ge T$ and
  $i \in \until{N}$.
  
  To prove the statement, we compare the state of the
  algorithm at an iteration $t \ge T_\delta$ and the quantities that would be
  computed at infinity, for all $i \in \until{N}$.
  To this end, let us denote by $(\algrho_i, \algx_i)$ the optimal solution of
  problem~\eqref{eq:alg_lexmin_MILP} with $\by_i^t = \by_i^\star$
  (we discard the $\xi_i$ part of the solution).
  As shown in
  Section~\ref{sec:mixed_integer_sol_computation}, in order to compute $\rho_i^t$,
  agents can equivalently solve
  \begin{align}
	\label{eq:ft_proof_rho_prob}
	\begin{split}
	  \min_{\rho_i, \bx_i} \: & \: \rho_i
	  \\
	  \subj \: 
	  & \: A_i \bx_i \leq  \by_i^t + \rho_i \1
	  \\
	  & \: \bx_i \in \Xi, \: \: \rho_i \ge 0.
	\end{split}
	\end{align}
  The pair $(\epsilon_i + \algrho_i, \algx_i)$ is feasible for
  problem~\eqref{eq:ft_proof_rho_prob} for all $t \ge T$.
  Indeed, it holds $\algx_i \in X_i$, $\epsilon_i + \algrho_i \ge 0$, and moreover
    $A_i\algx_i \le \by_i^\star + \algrho_i \1 
     \le \by_i^t + \epsilon_i \1 + \algrho_i \1$
  for all $t \ge T$.
  Being %
  $\rho_i^t $ the optimal cost
  of~\eqref{eq:ft_proof_rho_prob}, we have
  $\rho_i^t \le \epsilon_i + \algrho_i$ for all $t \ge T$.
  
  We now follow arguments similar to the proof of Theorem~\ref{thm:asymptotic_feasibility}.
  For all $i \in \until{N}$, let $\bz_i^\star$ denote the optimal solution
  of problem~\eqref{eq:LP_restricted}. For the sake of analysis, let us split
  the agent set $\until{N}$ as $\INT \cup  I_\feas \cup  I_\infeas$,
  where $\INT$ contains agents for which $\bz_i^\star \in \Xi$,
  $I_\feas$ contains agents for which $\bz_i^\star \notin \Xi$
  and $\algrho_i = 0$, and $I_\infeas$ contains agents
  for which $\bz_i^\star \notin \Xi$ and $\algrho_i > 0$.
  Using the same arguments of Theorem~\ref{thm:asymptotic_feasibility},
  for all $i \in \INT$ it holds $\algrho_i = 0$.
  
  Now, let us consider the agents $i \in \INT \cup I_\feas$. By construction, it holds
  \begin{align}
    A_i \bx_i^t
     \le \by_i^t + \rho_i^t
    \le \by_i^t + \epsilon_i,
    \:
    \forall \: i \in \INT \!\cup\! I_\feas \text{ and } t \ge T_\delta,
  \label{eq:ft_proof_feas}
  \end{align}
  where we also used $\rho_i^t \le \epsilon_i + \algrho_i$.
  As for the agents $i \in I_\infeas$, again by $\rho_i^t \le \epsilon_i + \algrho_i$, it holds
  $A_i \bx_i^t - \by_i^t \le \rho_i^t\1 \le \algrho_i\1 + \epsilon_i\1$ for all
  $t \ge T_\delta$, or equivalently
  $A_i \bx_i^t - \by_i^t -\epsilon_i\1 \le \algrho_i\1$,
  for all $t \ge T_\delta$.
  Moreover, note that, for $t \ge T_\delta$,
  it holds $A_i \bx_i^t - \by_i^t -\epsilon_i\1 \le A_i \bx_i^t - \by_i^t \le \max_{\bx_i \in \Xi} A_i \bx_i - \bL_i$,
  where $\max$ is intended component wise. Using the definition of $\sigma_i^\textsc{loc}$
  in Section~\ref{sec:restriction_approach} and rearranging the terms,
  we obtain
  \begin{align}
  \label{eq:ft_proof_infeas}
    A_i \bx_i^t \le \by_i^t + \sigma_i^\textsc{loc} + \epsilon_i\1,
    \hspace{0.5cm}
    \forall \: i \in I_\infeas \text{ and } t \ge T_\delta.
  \end{align}
  
  Finally, by using~\eqref{eq:ft_proof_feas} and~\eqref{eq:ft_proof_infeas}, we can write
  \begin{align*}
    \smallsum_{i=1}^N A_i \bx_i^t
    & 
    = \smallsum_{i=1}^N ( \by_i^t + \epsilon_i \1) + \smallsum_{i \in I_\infeas} \sigma_i^\textsc{loc}
    \\
    &
    \le
    b - \asyrestriction - \delta \1 + \smallsum_{i=1}^N \epsilon_i \1 + \asyrestriction,
    \hspace{1cm}
    \forall t \ge T,
  \end{align*}
  which follows since $\sum_{i=1}^N \by_i^t = \sum_{i=1}^N \by_i^0 = b - \asyrestriction - \delta\1$
  (cf. also Proposition~\ref{prop:DPD_convergence} (i)) and
  by~\eqref{eq:i_infeas_sum}.
  Since $\epsilon_i$ are arbitrary, choosing $\epsilon_i = \delta/N$ for all $i$
  implies $\sum_{i=1}^N A_i \bx_i^t \le b$ for all $t \ge T_\delta \triangleq T$, which
  concludes the proof.
\oprocend

\subsection{Proof of Theorem~\ref{thm:finite_time_performance}}

  Let $\{\by_1^t, \ldots, \by_N^t\}_{t \ge 0}$ denote the allocation vector
  sequence generated by Algorithm~\ref{alg:algorithm}.
  By following similar arguments as in the proof of Theorem~\ref{thm:finite_time_feasibility},
  we conclude that, for fixed $\epsilon_i > 0$, there exists a sufficiently large $T > 0$
  such that $\| \by_i^\star - \by_i^t \|_\infty \le \epsilon_i$ for all $t \ge T$ and
  $i \in \until{N}$.
  
  As done in \cite[Theorem 3.3]{vujanic2016decomposition},
  let us split the suboptimality bound as
  $\sum_{i=1}^N c_i^\top \bx_i^t \!-\! J^\MILP =
    \sum_{i=1}^N ( c_i^\top \bx_i^t - J_i^{\LP,t} )
     +
    \big(  \sum_{i=1}^N J_i^{\LP,t} -  J^{\LP,\ftrestriction} \big)
     + 
    \big( J^{\LP,\ftrestriction} \!-\! J^\LP \big)
    \!+\!
    \big( J^\LP \!-\! J^\MILP \big),$
  where $ J^{\LP,\ftrestriction}$ denotes the optimal cost of problem~\eqref{eq:LP_restricted}
  with $\restriction = \ftrestriction $.
  The first term $\sum_{i=1}^N ( c_i^\top \bx_i^t - J_i^{\LP,t} )$ can be explicitly computed.
  As for the last two terms, by following similar arguments as in \cite{vujanic2016decomposition},
  we conclude that $\big( J^{\LP,\ftrestriction} - J^\LP \big)
    +
    \big( J^\LP - J^\MILP \big) \le \Gamma \|\ftrestriction\|_\infty$.
  Let us analyze in detail the second term.

Notice that $J^{\LP,\ftrestriction}$ can be seen as the optimal cost of a perturbed version
of the problem having optimal cost $\sum_{i=1}^N J_i^{\LP,t}$, namely
the aggregate problem solved by the agents at iteration $t$, i.e.,
\begin{align}
\label{eq:finite_time_subopt_LP_t}
\begin{split}
  \min_{\substack{\bz_1,\ldots,\bz_N, \\ v_1, \ldots, v_N}} \:
  & \: \smallsum_{i =1}^N (c_i^\top \bz_i + M v_i)
  \\
  \subj \: 
  & \: A_i \bz_i \le \by_i^t + v_i \1, \hspace{2cm} \forall \: i,
  \\
  & \: \bz_i \in \conv{\Xi}, \:\: v_i \ge 0, \hspace{0.6cm} \forall \: i,
\end{split}
\end{align}
In particular, the constraints $A_i \bz_i \leq \by_i^t + v_i\1$ are perturbed
by $\by_i^\star - \by_i^t$ to obtain $A_i \bz_i \leq \by_i^\star + v_i\1$.
By applying perturbation theory~\cite{boyd2004convex}, we have
for all $t \ge T$
\begin{align*}
  \smallsum_{i=1}^N J_i^{\LP,t} \!-\! J^{\LP,\ftrestriction} \!
  \le
  \smallsum_{i=1}^N \| \by_i^\star - \by_i^t\|_\infty \| \bmu_i^{t}\|_1
  \le
  \smallsum_{i=1}^N \epsilon_i \| \bmu_i^{t}\|_1.
\eqoprocend
\end{align*}

\clearpage
\newpage

\twocolumn[
	\begin{@twocolumnfalse}
		{
		\centering
		\huge
		Supplement to the paper:\\ Distributed Primal Decomposition for Large-Scale MILPs\\[0.5em]
		}
		{
		\centering
		Andrea Camisa, \IEEEmembership{IEEE Student Member},
    Ivano Notarnicola, \IEEEmembership{IEEE Member},
    Giuseppe Notarstefano, \IEEEmembership{IEEE Member}\\[3em]
		}
	\end{@twocolumnfalse}
]

\subsection{Proof of Proposition~\ref{prop:DPD_convergence}}
\label{app:proof_DPD_convergence}

Let $M > 0$ be sufficiently large such that \cite[Lemma III.2]{notarnicola2017constraint} %
applies. To prove the result, we prove the equivalence of the distributed
algorithm~\eqref{eq:dpd_y_update} with the algorithm in \cite{notarnicola2017constraint}.
Therefore, let us first recall the algorithm in~\cite{notarnicola2017constraint}, which
reads as follows. Each agent $i$ maintains variables $\blambda_{ij}^t$, $j \in \nbrs_i$.
At each iteration, agents gather $\blambda_{ji}^t$ from $j \in \nbrs_i$ and compute
$\big( (\bz_i^t, v_i^t), \bmu_i^t \big)$ as a primal-dual optimal solution pair of
\begin{align} 
\label{eq:rsdd_z_LP}
\begin{split}
	\min_{\bz_i, v_i} \hspace{1cm} &\: c_i^\top \bz_{i} + M v_i
	\\
	\subj \hspace{0.1cm} 
	\: \bmu_i : & \: A_i \bz_i - \frac{b - \restriction}{N} + \!\sum_{j \in \nbrs_i} \!\!\big( \blambda_{ij}^t - \blambda_{ji}^t \big) \leq v_i\1
	\\
	& \: \bz_i \in \conv{\Xi}, \: \: v_i \ge 0.
\end{split}
\end{align}
Then, they gather $\bmu_j^t$ from $j \in \nbrs_i$ and update $\blambda_{ij}^t$
with
\begin{align} 
\label{eq:rsdd_lambda_update}
  \blambda_{ij}^{t+1} = \blambda_{ij}^t - \gamma^t \big( \bmu_i^t - \bmu_j^t \big)
  \hspace{0.7cm}
  \forall \: j \in \nbrs_i,
\end{align}
where $\gamma^t$ is the step size.

We now show that the update~\eqref{eq:rsdd_lambda_update}
is equivalent to~\eqref{eq:dpd_y_update} up to a change
of variable. To this end, let us define for all $t \ge 0$
\begin{align}
  \by_i^t \triangleq - \sum_{j \in \nbrs_i} \big( \blambda_{ij}^t - \blambda_{ji}^t \big) + \frac{b - \restriction}{N},
  \hspace{0.3cm}
  i \in \until{N}.
\label{eq:y_definition}
\end{align}
Then, it holds
\begin{align*}
  \sum_{i=1}^N \by_i^t
  &=
  \underbrace{\sum_{i=1}^N \sum_{j \in \nbrs_i} \big( \blambda_{ji}^t - \blambda_{ij}^t \big)}_{= \: \0}
    + \sum_{i=1}^N \frac{b - \restriction}{N}
  \\
  &=
  b - \restriction
  \hspace{4cm}
  \forall \: t \ge 0,
\end{align*}
which follows since the graph is undirected. This motivates the assumption
$\sum_{i=1}^N \by_i^0 = b - \restriction$ and proves (i).

To prove (ii), we simply note that the update of $\by_i^t$, as defined
in~\eqref{eq:y_definition}, reads
\begin{align*}
  \by_i^{t+1}
  &= \frac{b - \restriction}{N} + \sum_{j \in \nbrs_i} \big( \blambda_{ji}^{t+1} - \blambda_{ij}^{t+1} \big)
  \\
  &= \frac{b - \restriction}{N} + \sum_{j \in \nbrs_i} \big( \blambda_{ji}^t - \blambda_{ij}^t \big) +  2\gamma^t \sum_{j \in \nbrs_i} \big( \bmu_i^t - \bmu_j^t \big)
  \\
  &= \by_i^t + 2\gamma^t \sum_{j \in \nbrs_i} \big( \bmu_i^t - \bmu_j^t \big), 
  \hspace{0.7cm}
  i \in \until{N},
\end{align*}
where we point out that each $\bmu_i^t$ is a dual optimal solution of~\eqref{eq:rsdd_z_LP},
or, equivalently, a Lagrange multiplier of~\eqref{eq:dpd_z_LP}. %
Then, by defining the step-size sequence $\alpha^t \triangleq 2 \gamma^t$,
we see that the update~\eqref{eq:rsdd_lambda_update}
coincides with~\eqref{eq:dpd_y_update}.
As proven in \cite{notarnicola2017constraint}, the sequence
$\{(\blambda_{ij}^t)_{(i,j) \in \EE}\}_{t \ge 0}$ converges to an
optimal solution $(\blambda_{ij}^\star)_{(i,j) \in \EE}$ of a suitable reformulation
of~\eqref{eq:primal_decomp_master} %
(in terms of the variables $\blambda_{ij}$).
Therefore, the sequence $\{\by_1^t, \ldots, \by_N^t\}$ converges to some
$(\by_1^\star, \ldots, \by_N^\star)$, which is feasible for
problem~\eqref{eq:primal_decomp_master} %
(by (i)) and
cost-optimal (by optimality of $(\blambda_{ij}^\star)_{(i,j) \in \EE}$).
This concludes the proof of (ii).
As problem~\eqref{eq:rsdd_z_LP} is equivalent to~\eqref{eq:dpd_z_LP}, %
part (iii) follows by \cite[Theorem 2.4 (ii)]{notarnicola2017constraint}.
\oprocend

\subsection{Proof of Theorem~\ref{thm:asymptotic_performance_apriori}}

  Following the ideas in~\cite[Theorem 3.3]{vujanic2016decomposition}, let us split
  the bound as
	\begin{equation*}
	\begin{split}
	  & 
	  \sum_{i=1}^N c_i^\top \algx_i - J^\MILP =
	  \\
	  & \!
	  \sum_{i=1}^N (c_i^\top \algx_i \!-\! c_i^\top \bz_i^\star)
	  +
	  \Big( \sum_{i=1}^N c_i^\top \bz_i^\star \!-\! J^\LP \Big)
	  +
	  \big( J^\LP \!-\! J^\MILP \big) ,
	\end{split}
	\end{equation*}
	where $(\bz_1^\star, \ldots, \bz_N^\star)$ is the optimal solution of
	problem~\eqref{eq:LP_restricted} and $J^\LP$ denotes the optimal
	cost of problem~\eqref{eq:LP_restricted} with $\restriction = \0$.
	Next, we analyze each term independently.
	
  (i)
  $\sum_{i=1}^N (c_i^\top \algx_i \!-\! c_i^\top \bz_i^\star)$.
	By Proposition~\ref{prop:LP_integer_components}, there exists $\INT$, with
	$|\INT| \ge N-S$, such that $\bz_i^\star \in \Xi$ for all $i \in \INT$.
	Thus, for $i \in \INT$, it holds $\algx_i = \bz_i^\star$, implying
	$c_i^\top \algx_i \!-\! c_i^\top \bz_i^\star = 0$.
	Therefore, by defining $\NINT \triangleq \until{N} \setminus \INT$, the sum reduces to
	\begin{align}
	\label{eq:subopt_b_first}
	  \sum_{i=1}^N (c_i^\top \algx_i \!-\! c_i^\top \bz_i^\star)
	  =
	  \sum_{i \in \NINT} (c_i^\top \algx_i \!-\! c_i^\top \bz_i^\star),
	\end{align}
	with $|\NINT| \le S$. Since $c_i^\top \algx_i \le \max\limits_{\bx_i \in \Xi} c_i^\top \bx_i$ and
  $\min\limits_{\bx_i \in \Xi} c_i^\top \bx_i \le c_i^\top \bz_i^\star$, it follows that
  \begin{align*}
	  \sum_{i=1}^N (c_i^\top \algx_i \!-\! c_i^\top \bz_i^\star)
	  \le
	  \sum_{i \in \NINT} \gamma_i
	  \le
	  S \max_{i \in \until{N}} \gamma_i.
	\end{align*}

  (ii)
  $\sum_{i=1}^N c_i^\top \bz_i^\star - J^\LP$.
  By following similar arguments to~\cite[Theorem 3.3]{vujanic2016decomposition},
  one can show that
  \begin{align}
  \label{eq:subopt_b_second}
    \sum_{i=1}^N c_i^\top \bz_i^\star - J^\LP
    \le
    \frac{\|\asyrestriction\|_\infty}{\zeta}
    \sum_{i=1}^N ( c_i^\top \hz_i -  c_i^\top \bz_i^\star ),
  \end{align}
  where $(\hz_1, \ldots, \hz_N)$ is any Slater point  (cf. Assumption~\ref{ass:slater}).
  Since $c_i^\top \hz_i \le \max\limits_{\bx_i \in \Xi} c_i^\top \bx_i$ and
  $\min\limits_{\bx_i \in \Xi} c_i^\top \bx_i \le c_i^\top \bz_i^\star$, it follows that
  \begin{align*}
    \sum_{i=1}^N c_i^\top \bz_i^\star - J^\LP
    & \le
    \frac{\|\asyrestriction\|_\infty}{\zeta}
    \sum_{i=1}^N \gamma_i
    \\
    & \le
    \frac{N \|\asyrestriction\|_\infty}{\zeta}
    \max_{i \in \until{N}} \gamma_i.
  \end{align*}
  
  (iii) $J^\LP - J^\MILP$.
  Being $J^\LP$ the cost of~\eqref{eq:LP_restricted} with $\restriction = 0$,
  which is a relaxed version of~\eqref{eq:MILP},
  then $J^\LP - J^\MILP \le 0$.
  Combining the results above, the bound follows.
\oprocend

\subsection{Proof of Corollary~\ref{cor:asymptotic_performance}}

  It is sufficient to follow the same line of
  Theorem~\ref{thm:asymptotic_performance_apriori}, but stopping
  the proof of (i) at~\eqref{eq:subopt_b_first}
  and stopping the proof of (ii) at~\eqref{eq:subopt_b_second}.
\oprocend

\end{document}